\catcode`\@=11
\def\undefine#1{\let#1\undefined}
\def\newsymbol#1#2#3#4#5{\let\next@\relax
 \ifnum#2=\@ne\let\next@\msafam@\else
 \ifnum#2=\tw@\let\next@\msbfam@\fi\fi
 \mathchardef#1="#3\next@#4#5}
\def\mathhexbox@#1#2#3{\relax
 \ifmmode\mathpalette{}{\m@th\mathchar"#1#2#3}%
 \else\leavevmode\hbox{$\m@th\mathchar"#1#2#3$}\fi}
\def\hexnumber@#1{\ifcase#1 0\or 1\or 2\or 3\or 4\or 5\or 6\or 7\or 8\or
 9\or A\or B\or C\or D\or E\or F\fi}

\newdimen\ex@
\ex@.2326ex
\def\varinjlim{\mathop{\vtop{\ialign{##\crcr
 \hfil\rm lim\hfil\crcr\noalign{\nointerlineskip}\rightarrowfill\crcr
 \noalign{\nointerlineskip\kern-\ex@}\crcr}}}}
\def\varprojlim{\mathop{\vtop{\ialign{##\crcr
 \hfil\rm lim\hfil\crcr\noalign{\nointerlineskip}\leftarrowfill\crcr
 \noalign{\nointerlineskip\kern-\ex@}\crcr}}}}
\def\varliminf{\mathop{\underline{\vrule height\z@ depth.2exwidth\z@
 \hbox{\rm lim}}}}

\font\tenmsa=msam10
\font\sevenmsa=msam7
\font\fivemsa=msam5
\newfam\msafam
\textfont\msafam=\tenmsa
\scriptfont\msafam=\sevenmsa
\scriptscriptfont\msafam=\fivemsa
\edef\msafam@{\hexnumber@\msafam}
\mathchardef\dabar@"0\msafam@39
\def\dashrightarrow{\mathrel{\dabar@\dabar@\mathchar"0\msafam@4B}}
\def\dashleftarrow{\mathrel{\mathchar"0\msafam@4C\dabar@\dabar@}}

\font\tenmsb=msbm10
\font\sevenmsb=msbm7
\font\fivemsb=msbm5
\newfam\msbfam
\textfont\msbfam=\tenmsb
\scriptfont\msbfam=\sevenmsb
\scriptscriptfont\msbfam=\fivemsb
\edef\msbfam@{\hexnumber@\msbfam}
\def\Bbb#1{{\fam\msbfam\relax#1}}
\def\widehat#1{\setbox\z@\hbox{$\m@th#1$}%
 \ifdim\wd\z@>\tw@ em\mathaccent"0\msbfam@5B{#1}%
 \else\mathaccent"0362{#1}\fi}
\font\teneufm=eufm10
\font\seveneufm=eufm7
\font\fiveeufm=eufm5
\newfam\eufmfam
\textfont\eufmfam=\teneufm
\scriptfont\eufmfam=\seveneufm
\scriptscriptfont\eufmfam=\fiveeufm
\def\frak#1{{\fam\eufmfam\relax#1}}

\newsymbol\boxtimes 1202

\catcode`\@=12

\magnification=\magstep1
\font\title = cmr10 scaled \magstep2

\font\smallmath= cmmi7
\font\author = cmcsc10
\font\addr = cmti7
\font\byabs = cmr7

\parindent=1em
\baselineskip 15pt
\hsize=12.3 cm
\vsize=18.5 cm

\newcount\refcount
\newcount\seccount
\newcount\sscount
\newcount\eqcount
\newcount\boxcount
\newcount\testcount
\newcount\bibcount
\boxcount = 128
\seccount = -1

\def\proc#1#2{\advance\sscount by 1
	\medskip\goodbreak\noindent{\author #1}
	{\tenrm{\number\sscount}}:\ \ {\it #2}}
\def\nproc#1#2#3{\advance\sscount by 1\global
	\edef#1{#2\ \number\sscount}	
	\medskip\goodbreak\noindent{\author #2}
	{\tenrm{\number\sscount}}:\ \ {\it #3}}
\def\proof{\medskip\noindent{\it Proof:\ \ }}

\def\eql#1{\global\advance\eqcount by 1\global
	\edef#1{(\number\eqcount)}\leqno{#1}}
\def\ref#1#2{\advance\refcount by 1\global
	\edef#1{[\number\refcount]}\setbox\boxcount=
	\vbox{\item{[\number\refcount]}#2}\advance\boxcount by 1}
\def\biblio{{\frenchspacing
	\bigskip\goodbreak\centerline{\bf REFERENCES}\medskip
	\bibcount = 128\loop\ifnum\testcount < \refcount
	\goodbreak\advance\testcount by 1\box\bibcount
	\advance\bibcount by 1\vskip 4pt\repeat\medskip}}

\def\emph{\bf}
\def\colon{{:}\;}
\def\|{|\;}


\def\scirc{{\scriptstyle\circ}}

\def\Z{{\Bbb Z}}
\def\Q{{\Bbb Q}}
\def\F{{\Bbb F}}

\def\A{{\cal A}}

\def\Hom{{\rm Hom}}

\def\End{{\rm End}}
\def\ab{{\rm ab}}
\def\Gal{{\rm Gal}}

\def\GL{{\rm GL}}
\def\O{{\cal O}}
\def\E{{\cal E}}

\def\qed{\hfill\hbox{$\sqcup$\llap{$\sqcap$}}\medskip}

\def\p{{\frak p}}
\def\smp{{\scriptstyle\p}}

\def\Spec{{\rm Spec}\,}
\def\im{\mathop{\rm im}\nolimits}
\def\mod#1{(\mathop{\rm mod}\nolimits #1)}

\ref\AR{Ailon, N.; Runick, Z.: Torsion points on curves and common divisors of
	$a^k-1$ and $b^k-1$,
	{\tt arXiv:~math.NT/0202102 v2}, preprint, Feb.~28, 2002.}
\ref\BGK{Banaszak, G.; Gajda, W.; Kraso\'n, P.:
	A support problem for the intermediate Jacobians of $\ell$-adic representations,
	{\tt http://www.math.uiuc.edu/Algebraic-Number-Theory/0332}, preprint, Jan. 29, 2002.}
\ref\BLR{Bosch, S.; L\"utkebohmert, W.; Raynaud, M.:
	{\it N\'eron models}, Ergebnisse der Mathematik und ihrer Grenzgebiete, 21.
	Springer-Verlag, Berlin, 1990.}
\ref\CRS{Corrales-Rodrig\'a\~nez, C.; Schoof, R.:
	The support problem and its elliptic analogue.
	{\it J. Number Theory}  {\bf 64}  (1997),  no. 2, 276--290.}
\ref\Del{Deligne, P.: Conjectures de Tate et Shafarevitch, S\'eminaire Bourbaki 1983/84,
	$\rm n^o$~616, {\it Ast\'erisque} 121/122.}
\ref\KP{Khare, C.; Prasad, D.:  Reduction of homomorphisms mod $p$ and algebraicity,
	{\tt arXiv:~math.NT/0211004v1}, preprint, Nov.~1, 2002.}
\ref\Kow{Kowalski, E.: Some local-global applications of Kummer theory, preprint.}
\ref\Lang{Lang, S.: {\it Algebra}, second edition, Addison-Wesley, Menlo Park, CA, 1984.}
\ref\Ber{Serre, J-P.: Lettre \`a Daniel Bertrand du 8/6/1984, {\it Collected Papers Vol. IV.}}
\ref\Rib{Serre, J-P.: Lettre \`a Ken Ribet du 7/3/1986, {\it Collected Papers Vol. IV.}}
\ref\Wong{Wong, S.: Power residues on abelian varieties.  
	{\it Manuscripta Math.}  {\bf 102}  (2000),  129--138.}

\centerline{\title THE SUPPORT PROBLEM}
\medskip
\centerline{\title FOR ABELIAN VARIETIES}
\bigskip
\centerline{\byabs BY}
\medskip
\noindent{\author\hfill Michael Larsen\footnote*
{\tenrm Partially supported by NSF
Grant DMS-0100537}\hfill}
\medskip
\centerline{\addr Department of Mathematics, Indiana University}
\centerline{\addr Bloomington, IN 47405, USA}
\bigskip

\centerline{\byabs ABSTRACT}
\smallskip
{\byabs \narrower\narrower
\textfont1 = \smallmath
Let $A$ be an abelian variety over a number field $K$.  If $P$ and $Q$ are $K$-rational
points of $A$ such that the order of the (mod $\smp$) reduction of $Q$ divides the
order of the (mod $\smp$) reduction of $P$ for almost all prime ideals $\smp$, then there exists
a $K$-endomorphism $\phi$ of $A$ and a positive integer $k$ such that
$\phi(P)=kQ$.
\par}

\medskip
This note solves the support problem for abelian varieties over number fields, thus 
answering a question of C.~Corrales-Rodrig\'a\~nez and R.~Schoof \CRS. 
Recently,   G.~Banaszak, W.~Gajda, and P.~Kraso\'n \BGK\ 
and C.~Khare and D.~Prasad \KP\ have solved the problem for certain classes
of abelian varieties for which the images of the $\ell$-adic Galois representations can be
particularly well understood.  A number of other authors have also made progress
recently on closely related problems, including E.~Kowalski \Kow, S.~Wong \Wong,
and N.~Ailon and Z.~Rudnick \AR.

I would like to thank R.~Schoof for his helpful comments on 
earlier versions of this paper. 

The main result is as follows:

\nproc\Main{Theorem}{Let $K$ be a number field, $\O_K$ its ring of integers, and $\O$
the coordinate ring of an open subscheme of $\Spec\O_K$.  Let $\A$ be an abelian scheme 
over $\O$ and $P,Q\in\A(\O)$ arbitrary sections.  Suppose that for all $n\in\Z$
and all prime ideals $\p$ of $\O$, we have the implication
$$nP\equiv0\ \mod{\p}\Rightarrow nQ\equiv0\ \mod{\p}.\eql\hyp$$
Then there exist a positive integer $k$ and an endomorphism $\phi\in\End_{\O}(\A)$ such that
$$\phi(P)=kQ.\eql\conc$$
}

Note that as $\A$ is a N\'eron model of its generic fiber $A$ (\BLR~I~1.2/8), we have that
$\End_{\O}\A=\End_K A$.  We employ scheme notation only to make sense of the notion
of the reduction of a point of $A$ $\mod{\p}$.

It is clear that if $Q=\phi(P)$, the order of any reduction of $Q$ divides that of
the corresponding reduction of $P$.  One might ask whether the converse
is true or, in other words, whether one can strengthen \conc\ to ask
that $k=1$.  The following proposition shows that in general the answer is
negative:

\proc{Proposition}{There exist $\O$, $\A$, $P$, and $Q$ as above such that
\hyp\ holds
but $Q\not\in(\End_{\O}\A) P$.}

\proof Let $\O$ be a ring containing $1/2$.
Let $\E/\O$ be an elliptic curve with $\End_{\O}\E=\Z$
whose $2$-torsion is all $\O$-rational.
Let $T_1$ and $T_2$ denote distinct $2$-torsion points of $\E(\O)$, and let $R$
denote a point of infinite order in $\E(\O)$.  Let
$\A=\E^2$, $P=(R,R+T_1)$, and $Q=(R,R+T_2)$.  Then the reductions of $R$ and $R+T_1$ cannot
both have odd order (since $T_1$ has order exactly 2 in any reduction $\mod{\p}$), so
$P$ always has even order $\mod{\p}$.  Thus $nP\equiv0\ \mod{\p}$ implies
$2\mid n$ and therefore
$$nQ=(nR,nR) = nP\equiv0\ \mod{\p}.$$
On the other hand, $\End_{\O}\A=M_2(\Z)$, so no endomorphism of $\A$ sends $P$ to $Q$.\qed

Let $E=\End_\O \A$.
We begin by showing that \conc\ is implied by its $\mod{m}$ analogue for sufficiently
large $m$.

\nproc\Modular{Lemma}{Given $\O$, $\A$, and $E$ as above and $\O$-points $P$ and $Q$
of $\A$, either $P$ and $Q$ satisfy \conc\ or 
there exists $n$ such that for all $\phi\in E$ and all $m\ge n$,
$$\phi(P)-Q\not\in m \A(\O).$$
}

\proof
The lemma follows from the Mordell-Weil theorem and
the trivial fact that the image of $Q$ in the finitely generated abelian 
group $\A(\O)/E P$ is of finite order
if it is $m$-divisible for infinitely many values of $m$.\qed

Next, we prove two simple algebraic lemmas.

\nproc\Inter{Lemma}{Let $G$ be a group with normal subgroups $G_1$ and $G_2$ such that
$G/G_i$ is finite and abelian for $i=1,2$.
Let $\alpha$ be an automorphism of $G$ such that $\alpha(G_i)\subset G_i$ 
for $i=1,2$.  Suppose
$\alpha$ acts trivially on $G/G_1$ and as a scalar $m$ on $G/G_2$, where
$m-1$ is prime to $G/G_2$.  Then every coset of
$G_1$ meets every coset of $G_2$.}

\proof 
Applying Goursat's lemma (\Lang~I, Ex.) to the $\alpha$-equivariant map
$$\psi\colon G/(G_1\cap G_2)\to G/G_1\times G/G_2,$$
we find normal subgroups $H_1\supset G_1$ and $H_2\supset G_2$ of $G$ (automatically
$\alpha$-stable) such that
the image of $\psi$ is the pullback to $G/G_1\times G/G_2$ of the graph of an
$\alpha$-equivariant isomorphism $G/H_1\tilde\to G/H_2$.  By hypothesis, the two sides of
this isomorphism must be trivial, so $\psi$ is surjective, which proves the lemma.\qed

\nproc\Endo{Lemma}{
Let $M$ and $N$ be left modules of a ring $R$.  Suppose that $N$ is 
semisimple.  Let $\alpha,\beta\in\Hom_R(M,N)$ be such that $\ker\alpha\subset\ker\beta$.  
Then there exists $\gamma\in\End_R(N)$ such that $\beta=\gamma\scirc\alpha$.}

Let $M_\alpha=\ker\alpha$ and $M_\beta=\ker\beta$, so $M_\alpha\subset M_\beta$.  Let $N_\alpha\cong M/M_\alpha$
and $N_\beta\cong M/M_\beta$ denote the images of $\alpha$ and $\beta$.  Thus, $N_\beta$
is isomorphic to a quotient of $N_\alpha$.  As $N$ is semisimple, there is a projection map
$N\to N_\alpha$.  Composing this with the quotient map $N_\alpha\to N_\beta$ and the
inclusion $N_\beta\subset N$ we obtain the desired map $\gamma$.\qed

We remark that \Endo\ holds more generally for any abelian category.


%
%
%


We can now prove the main theorem.
Let $\rho_\ell\colon G_K\to\GL_{2g}(\Z_\ell)$ denote the $\ell$-adic Galois
representation given by the Tate module of $A$, and let $\bar\rho_\ell$ 
denote its $\mod{\ell}$ reduction.  
Let $G_n$ denote the Galois group of the field $K_n$ of $n$-torsion
points on $A$.  In particular, $G_\ell$ is the image of $\bar\rho_\ell$.
Let $M_\ell=\End_{\Z} \bigl(A[\ell](\bar K)\bigr)\cong M_{2g}(\F_\ell)$ 
denote the endomorphism ring of the
additive group of $\ell$-torsion points of $A$ over $\bar K$.  We
choose $\ell$ sufficiently large that it enjoys the following properties:
\smallskip

\item{(a)}The group of homotheties in $\rho_\ell(G_K)$ is of index $<\ell-1$ in $\Z_\ell^*$.
\item{(b)}The image $E_\ell$ of $E$ in $M_\ell$ and the subring of $M_\ell$ generated
	by $G_\ell$ are mutual centralizers.  In particular, both are semisimple algebras.
\item{(c)}If for some $\phi\in E$, one has $\phi(P)-Q\in\ell A(K)$, then $P$ and $Q$
	satisfy \conc.

\smallskip

Part (a) follows from a result of Serre (\Rib\ \S2).
Part (b) is a well-known folklore corollary of Faltings' proof of the Tate conjecture.
See \Ber~p.~24 for a statement.    We sketch a proof.  The endomorphism ring
$E$ acts on $H^1_{\rm sing}(A,\Z)$.
Let $E^*$ be the centralizer of $E$ in $\End_{\Z} H^1_{\rm sing}(A,\Z)$ and $E^{**}$ its
double centralizer.  As $E\otimes\Q$ is semisimple, $E^{**}\otimes \Q = E\otimes\Q$,
so $E$ is of finite index in $E^{**}$.  For $\ell$ sufficiently large, therefore, $E_\ell = E^{**}_\ell$.
The commutator map gives a homomorphism of abelian groups 
$M_{2g}(\Z)\to \Hom(E, M_{2g}(\Z))$ with kernel $E^*$.  
The sequence
$$0\to E^*\to M_{2g}(\Z)\to \Hom(E, M_{2g}(\Z))$$
remain exact after tensoring with $\F_\ell$ for $\ell\gg 0$.  
Therefore, the commutator of $E_\ell$ in $M_\ell$ is
$E^*_\ell$ for $\ell\gg 0$, and likewise the commutator of $E^*_\ell$ in $M_\ell$ is 
$E^{**}_\ell = E_\ell$ for $\ell\gg0$.  By the double commutant theorem, 
$E_\ell$ and $E^*_\ell$ are semisimple.
Now, \Del~2.7 asserts that for all $\ell\gg 0$, the centralizer of
$E\otimes\Z_\ell$ in the endomorphism ring  of the $\ell$-adic Tate module  
$T_\ell A = H^1_{\rm sing}(A,\Z)\otimes\Z_\ell$, 
is the image of $\Z_\ell[G_K]$, or in other words, 
$\im(\Z_\ell[G_K]\to\End(T_\ell A)) = E^*\otimes\Z_\ell,$
which implies (b).
Part (c) follows from \Modular.

The Kummer sequence for $A/K$ gives a natural $E_\ell$-equivariant
embedding
$$A(K)/\ell A(K)\hookrightarrow H^1(G_K, A[\ell](\bar K)) = H^1(G_K, A[\ell](K_\ell)).$$
By (a), the group $G_\ell$ contains a non-trivial subgroup $S_\ell$ which acts by scalar
multiplication on $A[\ell](K_\ell)$.
Since
$$A[\ell](K_\ell)^{S_\ell} = H^1(S_\ell,A[\ell](K_\ell))=0,$$
the inflation-restriction sequence
$$0\to H^1(G_\ell/S_\ell,A[\ell](K_\ell)^{S_\ell})\to H^1(G_\ell,A[\ell](K_\ell))
\to H^1(S_\ell,A[\ell](K_\ell))^{G_\ell/S_\ell}$$
implies $H^1(G_\ell,A[\ell](K_\ell)) = 0$.
The inflation-restriction sequence
$$0\to H^1(G_\ell,A[\ell](K_\ell))\to  H^1(G_K, A[\ell](K_\ell))
\to H^1(G_{K_\ell},A[\ell](K_\ell))^{G_\ell}$$
implies
$$A(K)/\ell A(K)\hookrightarrow \Hom(G_{K_\ell},A[\ell](K_\ell))^{G_\ell}
=\Hom_{\F_\ell[{G_\ell}]}(G_{K_\ell}^{\ab}\otimes\F_\ell,A[\ell](K_\ell))\eql\bdry$$
is injective.
For any $X\in A(K)$, we write $[X]$ for the class of the image of $X+\ell A(X)$ in
the right hand side of \bdry.

Let $V_\ell = G_{K_\ell}^{\ab}\otimes\F_\ell$.
Suppose that for all $\sigma\in V_\ell$, the condition 
$[Q](\sigma)=0$ implies $[P](\sigma)=0$.
Applying \Endo\ to the $\F_\ell[G_\ell]$-modules $M=V_\ell$
and $N = A[\ell](K_\ell)$, we obtain an $\F_\ell[G_\ell]$-module endomorphism
$\gamma$ of $N$ 
such that $\gamma\scirc[P]=[Q]$.  By (b), the endomorphism
$\gamma$ lies in the image of $E_\ell$, and lifting
it to an endomorphism $\phi\in E$, we conclude
$[\phi(P)-Q]=0$.  By \bdry, this means $\phi(P) -  Q\in\ell A(K)$,
and by (c), this implies \conc.

Therefore, we may assume that there exists $\sigma\in V_\ell$ with
$[Q](\sigma)=0$ and $[P](\sigma)\neq 0$.  
The pair $(P,Q)$ defines a $G_\ell$-equivariant map $V_\ell\to A[\ell](K_\ell)\times A[\ell](K_\ell)$.
The Galois action on $A[\ell^2](\bar K)$ defines a $G_\ell$-equivariant map
$V_\ell\to M_\ell$ since we have
$$\Gal(K_{\ell^2}/K_\ell)=\ker(G_{\ell^2}\to G_\ell) 
{\buildrel \rm log\over\hookrightarrow}
\ker\bigl(\End(A[\ell^2](\bar K))\to\End(A[\ell](\bar K)\bigr) = M_\ell.$$
By (a), there exists a non-trivial homothety in $G_\ell$.  It acts 
trivially on $M_\ell$ since the action of $G_\ell$ on $M_\ell$ is by conjugation,
and by definition, it acts as a non-trivial scalar on $A[\ell](K_\ell)\times A[\ell](K_\ell)$.  By \Inter,
the image of $V_\ell$ in $A[\ell](K_\ell)\times A[\ell](K_\ell)\times M_\ell$ is the product of its 
images in $A[\ell](K_\ell)\times A[\ell](K_\ell)$ and in $M_\ell$.  Applying (a) again, there exists 
$\sigma\in V_\ell$ such that $[P](\sigma)\neq 0$, $[Q](\sigma) = 0$, and $\sigma$
maps to a non-zero homothety in $M_\ell$.  

Let $K_{\ell^2}(\ell^{-1}P,\ell^{-1}Q)$ denote the extension of $K_\ell$ associated to
$$\ker V_\ell\to A[\ell](K_\ell)\times A[\ell](K_\ell)\times M_\ell;$$
thus  $K_{\ell^2}(\ell^{-1}P,\ell^{-1}Q)$ is the extension of $K$ generated by the
coordinates of all points $R\in A(\bar K)$ such that $\ell R\in\Z P + \Z Q + A[\ell](K_\ell)$.
By Cebotarev, we can fix a prime 
$\p$ of $\O$ which is unramified in $K_{\ell^2}(\ell^{-1}P,\ell^{-1}Q)$
and whose Frobenius conjugacy class in $\Gal(K_{\ell^2}(\ell^{-1}P,\ell^{-1}Q)/K)$ 
contains the image of $\sigma$ in $\Gal(K_{\ell^2}(\ell^{-1}P,\ell^{-1}Q)/K_\ell)$.  
Reducing $\mod{\p}$ we obtain
a finite field $\F_{\p}$ such that the $\ell$-primary part of $\A(\F_{\p})$ contains
$(\Z/\ell\Z)^{2g}$ (since the Frobenius at $\p$ fixes $K_\ell$) but has no
element of order $\ell^2$ (since the Frobenius at $\p$ acts as a non-trivial
homothety on $A[\ell^2](K_{\ell^2}(\ell^{-1}P,\ell^{-1}Q)) = A[\ell^2](\bar K)$.)
Moreover, the image of $P$ in $\A(\F_{\p})$ is not divisible by $\ell$,
but the image of $Q$ is.  
This means that the order of $P$ is divisible by $\ell$ but the order of $Q$ is prime 
to $\ell$, contrary to \hyp.
\qed

\proc{Corollary}{Let $K$ be a number field, $\O_K$ its ring of integers, and $\O$
the coordinate ring of an open subscheme of $\Spec\O_K$.  Let $\A_1,\A_2$ be abelian schemes 
over $\O$ and $P_i\in\A_i(\O)$ arbitrary sections.  Suppose that for all $n\in\Z$
and all prime ideals $\p$ of $\O$, we have the implication
$$nP_1\equiv0\ \mod{\p}\Rightarrow nP_2\equiv0\ \mod{\p}.$$
Then there exist a positive integer $k$ and an endomorphism 
$\psi\in\Hom_{\O}(\A_1,\A_2)$ such that
$$\psi(P_1) = kP_2.$$
}

\proof
Let $\A=\A_1\times\A_2$, $P=(P_1,0)$, $Q=(0,P_2)$.  Applying \Main, we conclude that
there exist a positive integer $k$ and an endomorphism
$$\phi\in\End_{\O}\A = \End_{\O}\A_1\times\End_{\O}\A_2\times\Hom_{\O}(\A_1,\A_2)\times
\Hom_{\O}(\A_2,\A_1)$$
such that $\phi(P)=kQ$.  
Letting $\psi$ denote the image of $\phi$ under projection to
$\Hom_{\O}(\A_1,\A_2)$, we obtain the corollary.
\qed

\biblio
\end